# Mappings of Conjugation of Quaternion Algebra

Aleks Kleyn

ABSTRACT. In the paper I considered mappings of conjugation of quaternion algebra. I proved the theorem that there is unique expansion of $R$-linear mapping of quaternion algebra relative to the given set of mappings of conjugation.

## Contents

## 1. Preface

The theorems [1]-2.5.1, [3]-1.5.1 state that $R$-linear mapping of complex field has form

(1.1) $$f = a + b \circ I$$
(1.2) $$(a + b \circ I) \circ z = az + b\overline{z}$$

where
$$I \circ z = \overline{z}$$

is conjugation of complex field. The identity map is linear homomorphism. Conjugation is antilinear homomorphism.

In the paper [2], I discovered that there exist nontrivial linear automorphisms of quaternion algebra. It inspired me to assume that it is possible to expand $R$-linear mapping into a sum of linear and antilinear automorphisms. Such an expansion exists (the theorem [2]-12.1 ). However such an expansion of linear mapping does not allow us to select in the set of linear mappings a class of mappings that would be similar to the linear homomorphism (the example [2]-12.3 ).

Aleks_Kleyn@MailAPS.org.
http://sites.google.com/site/AleksKleyn/.
http://arxiv.org/a/kleyn_a_1.
http://AleksKleyn.blogspot.com/.

The dimension of complex field equals 2, however the dimension of quaternion algebra equals 4. So a task to find a satisfactory expansion of the linear mapping relative to the basis consisting of linear and antilinear homomorphisms is not easy task.

The mapping $I$ in the equation (1.1) is conjugation. However we cannot find expansion of linear mapping into sum of identity mapping and conjugation, because the set of these mappings is too small to construct the basis.

Conjugation is a mapping of symmetry. Under conjugation, the coefficients of the $i$, $j$, $k$ change sign. In this paper, I want to explore how it will look the expansion of linear mapping, if I select a basis of mappings

$$(\overline{E}, \overline{I}, \overline{J}, \overline{K})$$

where each of these mappings change sign of only one coefficient.

## 2. Conventions

**Convention 2.1.** *Let A be free finite dimensional algebra. Considering expansion of element of algebra A relative basis $\overline{\overline{e}}$ we use the same root letter to denote this element and its coordinates. However we do not use vector notation in algebra. In expression $a^2$, it is not clear whether this is component of expansion of element a relative basis, or this is operation $a^2 = aa$. To make text clearer we use separate color for index of element of algebra. For instance,*

$$a = a^i \overline{e}_i$$

□

**Convention 2.2.** *If free finite dimensional algebra has unit, then we identify the vector of basis $\overline{e}_0$ with unit of algebra.* □

Without a doubt, the reader may have questions, comments, objections. I will appreciate any response.

## 3. Mapping $\overline{E}$

Consider mapping

$$\overline{E} : H \to H \quad \overline{E} \circ x = x$$

$$E^i_j = \delta^i_j$$

**Theorem 3.1.** *We can identify the mapping*

$$a \circ \overline{E} : H \to H \quad a \in H$$

*and matrix*

(3.1) $$E_{l \cdot a} = \begin{pmatrix} a^0 & -a^1 & -a^2 & -a^3 \\ a^1 & a^0 & -a^3 & a^2 \\ a^2 & a^3 & a^0 & -a^1 \\ a^3 & -a^2 & a^1 & a^0 \end{pmatrix}$$

*Proof.* The product of quaternions
$$a = a^0 + a^1 i + a^2 j + a^3 k$$
and
$$x = x^0 + x^1 i + x^2 j + x^3 k$$
has form
$$ax = a^0 x^0 - a^1 x^1 - a^2 x^2 - a^3 x^3 + (a^0 x^1 + a^1 x^0 + a^2 x^3 - a^3 x^2)i$$
$$+ (a^0 x^2 + a^2 x^0 + a^3 x^1 - a^1 x^3)j + (a^0 x^3 + a^3 x^0 + a^1 x^2 - a^2 x^1)k$$
Therefore, function $f_a(x) = ax$ has Jacobian matrix (3.1). $\square$

**Theorem 3.2.** *We can identify the mapping*
$$a \star \overline{E} : H \to H \quad a \in H$$
*and matrix*

(3.2)
$$E_{r \cdot a} = \begin{pmatrix} a^0 & -a^1 & -a^2 & -a^3 \\ a^1 & a^0 & a^3 & -a^2 \\ a^2 & -a^3 & a^0 & a^1 \\ a^3 & a^2 & -a^1 & a^0 \end{pmatrix}$$

*Proof.* The product of quaternions
$$x = x^0 + x^1 i + x^2 j + x^3 k$$
and
$$a = a^0 + a^1 i + a^2 j + a^3 k$$
has form
$$xa = x^0 a^0 - x^1 a^1 - x^2 a^2 - x^3 a^3 + (x^0 a^1 + x^1 a^0 + x^2 a^3 - x^3 a^2)i$$
$$+ (x^0 a^2 + x^2 a^0 + x^3 a^1 - x^1 a^3)j + (x^0 a^3 + x^3 a^0 + x^1 a^2 - x^2 a^1)k$$
Therefore, function $f_a(x) = xa$ has Jacobian matrix (3.2). $\square$

## 4. Mapping $\overline{I}$

Consider mapping
$$\overline{I} : H \to H \quad \overline{I} \circ x = x^0 - x^1 i + x^2 j + x^3 k$$

$$I^0_0 = 1 \quad I^1_1 = -1 \quad I^2_2 = 1 \quad I^3_3 = 1$$

$$I = \begin{pmatrix} 1 & 0 & 0 & 0 \\ 0 & -1 & 0 & 0 \\ 0 & 0 & 1 & 0 \\ 0 & 0 & 0 & 1 \end{pmatrix}$$

**Theorem 4.1.** *We can identify the mapping*

$$a \circ \overline{I} : H \to H \quad a \in H$$

*and matrix*

(4.1) $$I_{l \cdot a} = \begin{pmatrix} a^0 & a^1 & -a^2 & -a^3 \\ a^1 & -a^0 & -a^3 & a^2 \\ a^2 & -a^3 & a^0 & -a^1 \\ a^3 & a^2 & a^1 & a^0 \end{pmatrix}$$

(4.2) $$I_{l \cdot a} = E_{l \cdot a \circ} {}^\circ I$$

*Proof.* The product of quaternions

$$a = a^0 + a^1 i + a^2 j + a^3 k$$

and

$$I \circ x = x^0 - x^1 i + x^2 j + x^3 k$$

has form

$$a \circ I \circ x = a^0 x^0 + a^1 x^1 - a^2 x^2 - a^3 x^3 + (-a^0 x^1 + a^1 x^0 + a^2 x^3 - a^3 x^2)i$$
$$+ (a^0 x^2 + a^2 x^0 - a^3 x^1 - a^1 x^3)j + (a^0 x^3 + a^3 x^0 + a^1 x^2 + a^2 x^1)k$$

Therefore, function $f_a \circ x = a \circ I \circ x$ has Jacobian matrix (4.1). The equation (4.2) follows from the chain of equations

$$E_{l \cdot a \circ} {}^\circ I = \begin{pmatrix} a^0 & -a^1 & -a^2 & -a^3 \\ a^1 & a^0 & -a^3 & a^2 \\ a^2 & a^3 & a^0 & -a^1 \\ a^3 & -a^2 & a^1 & a^0 \end{pmatrix} \circ \begin{pmatrix} 1 & 0 & 0 & 0 \\ 0 & -1 & 0 & 0 \\ 0 & 0 & 1 & 0 \\ 0 & 0 & 0 & 1 \end{pmatrix}$$

$$= \begin{pmatrix} a^0 & a^1 & -a^2 & -a^3 \\ a^1 & -a^0 & -a^3 & a^2 \\ a^2 & -a^3 & a^0 & -a^1 \\ a^3 & a^2 & a^1 & a^0 \end{pmatrix} = I_{l \cdot a}$$

□

**Theorem 4.2.** *We can identify the mapping*

$$a \star \overline{I} : H \to H \quad a \in H$$

*and matrix*

(4.3) $$I_{r \cdot a} = \begin{pmatrix} a^0 & a^1 & -a^2 & -a^3 \\ a^1 & -a^0 & a^3 & -a^2 \\ a^2 & a^3 & a^0 & a^1 \\ a^3 & -a^2 & -a^1 & a^0 \end{pmatrix}$$

(4.4) $$I_{r \cdot a} = E_{r \cdot a} \circ I$$

*Proof.* The product of quaternions
$$I \circ x = x^0 - x^1 i + x^2 j + x^3 k$$
and
$$a = a^0 + a^1 i + a^2 j + a^3 k$$
has form
$$(I \circ x)a = x^0 a^0 + x^1 a^1 - x^2 a^2 - x^3 a^3 + (x^0 a^1 - x^1 a^0 + x^2 a^3 - x^3 a^2)i$$
$$+ (x^0 a^2 + x^2 a^0 + x^3 a^1 + x^1 a^3)j + (x^0 a^3 + x^3 a^0 - x^1 a^2 - x^2 a^1)k$$

Therefore, function $f_a \circ x = a \star I \circ x$ has Jacobian matrix (4.3). The equation (4.4) follows from the chain of equations

$$E_{r \cdot a} \circ I = \begin{pmatrix} a^0 & -a^1 & -a^2 & -a^3 \\ a^1 & a^0 & a^3 & -a^2 \\ a^2 & -a^3 & a^0 & a^1 \\ a^3 & a^2 & -a^1 & a^0 \end{pmatrix} \circ \begin{pmatrix} 1 & 0 & 0 & 0 \\ 0 & -1 & 0 & 0 \\ 0 & 0 & 1 & 0 \\ 0 & 0 & 0 & 1 \end{pmatrix}$$

$$= \begin{pmatrix} a^0 & a^1 & -a^2 & -a^3 \\ a^1 & -a^0 & a^3 & -a^2 \\ a^2 & a^3 & a^0 & a^1 \\ a^3 & -a^2 & -a^1 & a^0 \end{pmatrix} = I_{r \cdot a}$$

□

## 5. Mapping $\overline{J}$

Consider mapping

$$\overline{J} : H \to H \quad \overline{J} \circ x = x^0 + x^1 i - x^2 j + x^3 k$$

$$J_0^0 = 1 \quad J_1^1 = 1 \quad J_2^2 = -1 \quad J_3^3 = 1$$

$$J = \begin{pmatrix} 1 & 0 & 0 & 0 \\ 0 & 1 & 0 & 0 \\ 0 & 0 & -1 & 0 \\ 0 & 0 & 0 & 1 \end{pmatrix}$$

**Theorem 5.1.** *We can identify the mapping*

$$a \circ \overline{J} : H \to H \quad a \in H$$

*and matrix*

(5.1) $$J_{l \cdot a} = \begin{pmatrix} a^0 & -a^1 & a^2 & -a^3 \\ a^1 & a^0 & a^3 & a^2 \\ a^2 & a^3 & -a^0 & -a^1 \\ a^3 & -a^2 & -a^1 & a^0 \end{pmatrix}$$

(5.2) $$J_{l \cdot a} = E_{l \cdot a \circ} \circ J$$

*Proof.* The product of quaternions
$$a = a^0 + a^1 i + a^2 j + a^3 k$$
and
$$J \circ x = x^0 + x^1 i - x^2 j + x^3 k$$
has form
$$a \circ J \circ x = a^0 x^0 - a^1 x^1 + a^2 x^2 - a^3 x^3 + (a^0 x^1 + a^1 x^0 + a^2 x^3 + a^3 x^2) i$$
$$+ (-a^0 x^2 + a^2 x^0 + a^3 x^1 - a^1 x^3) j + (a^0 x^3 + a^3 x^0 - a^1 x^2 - a^2 x^1) k$$

Therefore, function $f_a \circ x = a \circ J \circ x$ has Jacobian matrix (5.1). The equation (5.2) follows from the chain of equations

$$E_{l \cdot a \circ} \circ J = \begin{pmatrix} a^0 & -a^1 & -a^2 & -a^3 \\ a^1 & a^0 & -a^3 & a^2 \\ a^2 & a^3 & a^0 & -a^1 \\ a^3 & -a^2 & a^1 & a^0 \end{pmatrix} \circ \begin{pmatrix} 1 & 0 & 0 & 0 \\ 0 & 1 & 0 & 0 \\ 0 & 0 & -1 & 0 \\ 0 & 0 & 0 & 1 \end{pmatrix}$$

$$= \begin{pmatrix} a^0 & -a^1 & a^2 & -a^3 \\ a^1 & a^0 & a^3 & a^2 \\ a^2 & a^3 & -a^0 & -a^1 \\ a^3 & -a^2 & -a^1 & a^0 \end{pmatrix} = J_{l \cdot a}$$

□

**Theorem 5.2.** *We can identify the mapping*
$$a \star \overline{J} : H \to H \quad a \in H$$
*and matrix*

(5.3) $$J_{r \cdot a} = \begin{pmatrix} a^0 & -a^1 & a^2 & -a^3 \\ a^1 & a^0 & -a^3 & -a^2 \\ a^2 & -a^3 & -a^0 & a^1 \\ a^3 & a^2 & a^1 & a^0 \end{pmatrix}$$

(5.4) $$J_{r \cdot a} = E_{r \cdot a \circ} \circ J$$

*Proof.* The product of quaternions
$$J \circ x = x^0 + x^1 i - x^2 j + x^3 k$$
and
$$a = a^0 + a^1 i + a^2 j + a^3 k$$
has form
$$(J \circ x)a = x^0 a^0 - x^1 a^1 + x^2 a^2 - x^3 a^3 + (x^0 a^1 + x^1 a^0 - x^2 a^3 - x^3 a^2)i$$
$$+ (x^0 a^2 - x^2 a^0 + x^3 a^1 - x^1 a^3)j + (x^0 a^3 + x^3 a^0 + x^1 a^2 + x^2 a^1)k$$

Therefore, function $f_a \circ x = a \star J \circ x$ has Jacobian matrix (5.3). The equation (5.4) follows from the chain of equations

$$E_{r \cdot a} \circ I = \begin{pmatrix} a^0 & -a^1 & -a^2 & -a^3 \\ a^1 & a^0 & a^3 & -a^2 \\ a^2 & -a^3 & a^0 & a^1 \\ a^3 & a^2 & -a^1 & a^0 \end{pmatrix} \circ \begin{pmatrix} 1 & 0 & 0 & 0 \\ 0 & -1 & 0 & 0 \\ 0 & 0 & 1 & 0 \\ 0 & 0 & 0 & 1 \end{pmatrix}$$

$$= \begin{pmatrix} a^0 & a^1 & -a^2 & -a^3 \\ a^1 & -a^0 & a^3 & -a^2 \\ a^2 & a^3 & a^0 & a^1 \\ a^3 & -a^2 & -a^1 & a^0 \end{pmatrix} = I_{r \cdot a}$$

□

## 6. Mapping $\overline{K}$

Consider mapping

$$\overline{K} : H \to H \quad \overline{K} \circ x = x^0 + x^1 i + x^2 j - x^3 k$$

$$K_0^0 = 1 \quad K_1^1 = 1 \quad K_2^2 = 1 \quad K_3^3 = -1$$

$$K = \begin{pmatrix} 1 & 0 & 0 & 0 \\ 0 & 1 & 0 & 0 \\ 0 & 0 & 1 & 0 \\ 0 & 0 & 0 & -1 \end{pmatrix}$$

**Theorem 6.1.** *We can identify the mapping*

$$a \circ \overline{K} : H \to H \quad a \in H$$

*and matrix*

(6.1)
$$K_{l \cdot a} = \begin{pmatrix} a^0 & -a^1 & -a^2 & a^3 \\ a^1 & a^0 & -a^3 & -a^2 \\ a^2 & a^3 & a^0 & a^1 \\ a^3 & -a^2 & a^1 & -a^0 \end{pmatrix}$$

(6.2) $$K_{l \cdot a} = E_{l \cdot a \circ} {}^{\circ} K$$

*Proof.* The product of quaternions
$$a = a^0 + a^1 i + a^2 j + a^3 k$$
and
$$K \circ x = x^0 + x^1 i + x^2 j - x^3 k$$
has form
$$a \circ K \circ x = a^0 x^0 - a^1 x^1 - a^2 x^2 + a^3 x^3 + (a^0 x^1 + a^1 x^0 - a^2 x^3 - a^3 x^2) i$$
$$+ (a^0 x^2 + a^2 x^0 + a^3 x^1 + a^1 x^3) j + (-a^0 x^3 + a^3 x^0 + a^1 x^2 - a^2 x^1) k$$

Therefore, function $f_a \circ x = a \circ K \circ x$ has Jacobian matrix (6.1). The equation (6.2) follows from the chain of equations

$$E_{l \cdot a \circ} {}^{\circ} K = \begin{pmatrix} a^0 & -a^1 & -a^2 & -a^3 \\ a^1 & a^0 & -a^3 & a^2 \\ a^2 & a^3 & a^0 & -a^1 \\ a^3 & -a^2 & a^1 & a^0 \end{pmatrix} {}^{\circ}_{\circ} \begin{pmatrix} 1 & 0 & 0 & 0 \\ 0 & 1 & 0 & 0 \\ 0 & 0 & 1 & 0 \\ 0 & 0 & 0 & -1 \end{pmatrix}$$

$$= \begin{pmatrix} a^0 & -a^1 & -a^2 & a^3 \\ a^1 & a^0 & -a^3 & -a^2 \\ a^2 & a^3 & a^0 & a^1 \\ a^3 & -a^2 & a^1 & -a^0 \end{pmatrix} = K_{l \cdot a}$$

□

**Theorem 6.2.** *We can identify the mapping*
$$a \star \overline{K} : H \to H \quad a \in H$$
*and matrix*

(6.3) $$K_{r \cdot a} = \begin{pmatrix} a^0 & -a^1 & -a^2 & a^3 \\ a^1 & a^0 & a^3 & a^2 \\ a^2 & -a^3 & a^0 & -a^1 \\ a^3 & a^2 & -a^1 & -a^0 \end{pmatrix}$$

(6.4) $$K_{r \cdot a} = E_{r \cdot a \circ} {}^{\circ} K$$

*Proof.* The product of quaternions
$$K \circ x = x^0 + x^1 i + x^2 j - x^3 k$$
and
$$a = a^0 + a^1 i + a^2 j + a^3 k$$
has form
$$(K \circ x) a = x^0 a^0 - x^1 a^1 - x^2 a^2 + x^3 a^3 + (x^0 a^1 + x^1 a^0 + x^2 a^3 + x^3 a^2) i$$
$$+ (x^0 a^2 + x^2 a^0 - x^3 a^1 - x^1 a^3) j + (x^0 a^3 - x^3 a^0 + x^1 a^2 - x^2 a^1) k$$

Therefore, function $f_a \circ x = a \star K \circ x$ has Jacobian matrix (6.3). The equation (6.4) follows from the chain of equations

$$E_{r \cdot a} \circ I = \begin{pmatrix} a^0 & -a^1 & -a^2 & -a^3 \\ a^1 & a^0 & a^3 & -a^2 \\ a^2 & -a^3 & a^0 & a^1 \\ a^3 & a^2 & -a^1 & a^0 \end{pmatrix} \circ \begin{pmatrix} 1 & 0 & 0 & 0 \\ 0 & -1 & 0 & 0 \\ 0 & 0 & 1 & 0 \\ 0 & 0 & 0 & 1 \end{pmatrix}$$

$$= \begin{pmatrix} a^0 & a^1 & -a^2 & -a^3 \\ a^1 & -a^0 & a^3 & -a^2 \\ a^2 & a^3 & a^0 & a^1 \\ a^3 & -a^2 & -a^1 & a^0 \end{pmatrix} = I_{r \cdot a}$$

$\square$

## 7. Linear Mapping

**Theorem 7.1.** *Linear mapping of quaternion algebra*

$$\overline{f} : H \to H$$

*has unique expansion*

(7.1) $$\overline{f} = a_0 \circ \overline{E} + a_1 \circ \overline{I} + a_2 \circ \overline{J} + a_3 \circ \overline{K}$$

*In this case*

(7.2) $$a_0^0 = \frac{1}{2}(-f_0^0 + f_1^1 + f_2^2 + f_3^3)$$

(7.3) $$a_1^0 = \frac{1}{2}(f_0^0 - f_1^1)$$

(7.4) $$a_2^0 = \frac{1}{2}(f_0^0 - f_2^2)$$

(7.5) $$a_3^0 = \frac{1}{2}(f_0^0 - f_3^3)$$

(7.6) $$a_0^1 = \frac{1}{2}(-f_0^1 - f_1^0 + f_2^3 - f_3^2)$$

(7.7) $$a_1^1 = \frac{1}{2}(f_0^1 + f_1^0)$$

(7.8) $$a_2^1 = \frac{1}{2}(f_0^1 - f_2^3)$$

(7.9) $$a_3^1 = \frac{1}{2}(f_0^1 + f_3^2)$$

$$(7.10) \quad a_0^2 = \frac{1}{2}(-f_0^2 - f_1^3 - f_2^0 + f_3^1)$$

$$(7.11) \quad a_1^2 = \frac{1}{2}(f_0^2 + f_1^3)$$

$$(7.12) \quad a_2^2 = \frac{1}{2}(f_0^2 + f_2^0)$$

$$(7.13) \quad a_3^2 = \frac{1}{2}(f_0^2 - f_3^1)$$

$$(7.14) \quad a_0^3 = \frac{1}{2}(-f_0^3 + f_1^2 - f_2^1 - f_3^0)$$

$$(7.15) \quad a_1^3 = \frac{1}{2}(f_0^3 - f_1^2)$$

$$(7.16) \quad a_2^3 = \frac{1}{2}(f_0^3 + f_2^1)$$

$$(7.17) \quad a_3^3 = \frac{1}{2}(f_0^3 + f_3^0)$$

*Proof.* Linear mapping (7.1) has matrix

$$(7.18) \quad \begin{pmatrix} a_0^0 & -a_0^1 & -a_0^2 & -a_0^3 \\ a_0^1 & a_0^0 & -a_0^3 & a_0^2 \\ a_0^2 & a_0^3 & a_0^0 & -a_0^1 \\ a_0^3 & -a_0^2 & a_0^1 & a_0^0 \end{pmatrix} + \begin{pmatrix} a_1^0 & a_1^1 & -a_1^2 & -a_1^3 \\ a_1^1 & -a_1^0 & -a_1^3 & a_1^2 \\ a_1^2 & -a_1^3 & a_1^0 & -a_1^1 \\ a_1^3 & a_1^2 & a_1^1 & a_1^0 \end{pmatrix}$$
$$+ \begin{pmatrix} a_2^0 & -a_2^1 & a_2^2 & -a_2^3 \\ a_2^1 & a_2^0 & a_2^3 & a_2^2 \\ a_2^2 & a_2^3 & -a_2^0 & -a_2^1 \\ a_2^3 & -a_2^2 & -a_2^1 & a_2^0 \end{pmatrix} + \begin{pmatrix} a_3^0 & -a_3^1 & -a_3^2 & a_3^3 \\ a_3^1 & a_3^0 & -a_3^3 & -a_3^2 \\ a_3^2 & a_3^3 & a_3^0 & a_3^1 \\ a_3^3 & -a_3^2 & a_3^1 & -a_3^0 \end{pmatrix}$$

From a comparison of matrix of the mapping $\overline{f}$ and matrix (7.18), we get the system of linear equations

$$(7.19) \quad f_0^0 = a_0^0 + a_1^0 + a_2^0 + a_3^0 \qquad (7.27) \quad f_0^2 = a_0^2 + a_1^2 + a_2^2 + a_3^2$$

$$(7.20) \quad f_1^0 = -a_0^1 + a_1^1 - a_2^1 - a_3^1 \qquad (7.28) \quad f_1^2 = a_0^3 - a_1^3 + a_2^3 + a_3^3$$

$$(7.21) \quad f_2^0 = -a_0^2 - a_1^2 + a_2^2 - a_3^2 \qquad (7.29) \quad f_2^2 = a_0^0 + a_1^0 - a_2^0 + a_3^0$$

$$(7.22) \quad f_3^0 = -a_0^3 - a_1^3 - a_2^3 + a_3^3 \qquad (7.30) \quad f_3^2 = -a_0^1 - a_1^1 - a_2^1 + a_3^1$$

$$(7.23) \quad f_0^1 = a_0^1 + a_1^1 + a_2^1 + a_3^1 \qquad (7.31) \quad f_0^3 = a_0^3 + a_1^3 + a_2^3 + a_3^3$$

$$(7.24) \quad f_1^1 = a_0^0 - a_1^0 + a_2^0 + a_3^0 \qquad (7.32) \quad f_1^3 = -a_0^2 + a_1^2 - a_2^2 - a_3^2$$

$$(7.25) \quad f_2^1 = -a_0^3 - a_1^3 + a_2^3 - a_3^3 \qquad (7.33) \quad f_2^3 = a_0^1 + a_1^1 - a_2^1 + a_3^1$$

$$(7.26) \quad f_3^1 = a_0^2 + a_1^2 + a_2^2 - a_3^2 \qquad (7.34) \quad f_3^3 = a_0^0 + a_1^0 + a_2^0 - a_3^0$$

The equation (7.3) follows from equations (7.19), (7.24). The equation (7.4) follows from equations (7.19), (7.29). The equation (7.5) follows from equations (7.19), (7.34). The equation (7.2) follows from equations (7.3), (7.4), (7.5), (7.19).

The equation (7.7) follows from equations (7.23), (7.20). The equation (7.8) follows from equations (7.23), (7.33). The equation (7.9) follows from equations (7.23), (7.30). The equation (7.6) follows from equations (7.7), (7.8), (7.9), (7.23).

The equation (7.11) follows from equations (7.27), (7.32). The equation (7.12) follows from equations (7.27), (7.21). The equation (7.13) follows from equations (7.27), (7.26). The equation (7.10) follows from equations (7.11), (7.12), (7.13), (7.27).

The equation (7.15) follows from equations (7.31), (7.28). The equation (7.16) follows from equations (7.31), (7.25). The equation (7.17) follows from equations (7.31), (7.22). The equation (7.14) follows from equations (7.15), (7.16), (7.17), (7.31). $\square$

It follows from the theorem [3]-2.6.4 (you can see also the theorem [1]-3.7.4) that the set of linear endomorphisms $\mathcal{L}(H;H)$ of quaternion algebra $H$ is isomorphic to tensor product $H \otimes H$. The theorem 7.1 states that we can consider the module $\mathcal{L}(H;H)$ as $H\star$-vector space with basis

$$(\overline{E}, \overline{I}, \overline{J}, \overline{K})$$

**Example 7.2.** According to the theorem 3.2, mapping

$$a \star \overline{E} : H \to H \quad a \in H$$

has matrix

$$\begin{pmatrix} a^0 & -a^1 & -a^2 & -a^3 \\ a^1 & a^0 & a^3 & -a^2 \\ a^2 & -a^3 & a^0 & a^1 \\ a^3 & a^2 & -a^1 & a^0 \end{pmatrix}$$

According to the theorem 7.1,

$$a_0^0 = \frac{1}{2}(-a^0 + a^0 + a^0 + a^0) \qquad a_1^0 = \frac{1}{2}(a^0 - a^0)$$

$$a_0^1 = \frac{1}{2}(-a^1 - (-a^1) + (-a^1) - a^1) \qquad a_1^1 = \frac{1}{2}(a^1 + (-a^1))$$

$$a_0^2 = \frac{1}{2}(-a^2 - a^2 - (-a^2) + (-a^2)) \qquad a_1^2 = \frac{1}{2}(a^2 + a^2)$$

$$a_0^3 = \frac{1}{2}(-a^3 + (-a^3) - a^3 - (-a^3)) \qquad a_1^3 = \frac{1}{2}(a^3 - (-a^3))$$

$$a_2^0 = \frac{1}{2}(a^0 - a^0) \qquad a_3^0 = \frac{1}{2}(a^0 - a^0)$$

$$a_2^1 = \frac{1}{2}(a^1 - (-a^1)) \qquad a_3^1 = \frac{1}{2}(a^1 + a^1)$$

$$a_2^2 = \frac{1}{2}(a^2 + (-a^2)) \qquad a_3^2 = \frac{1}{2}(a^2 - (-a^2))$$

$$a_2^3 = \frac{1}{2}(a^3 + a^3) \qquad a_3^3 = \frac{1}{2}(a^3 + (-a^3))$$

Therefore

$$a_0^0 = a^0 \qquad a_1^0 = 0 \qquad a_2^0 = 0 \qquad a_3^0 = 0$$
$$a_0^1 = -a^1 \qquad a_1^1 = 0 \qquad a_2^1 = a^1 \qquad a_3^1 = a^1$$
$$a_0^2 = -a^2 \qquad a_1^2 = a^2 \qquad a_2^2 = 0 \qquad a_3^2 = a^2$$
$$a_0^3 = -a^3 \qquad a_1^3 = a^3 \qquad a_2^3 = a^3 \qquad a_3^3 = 0$$

□

**Example 7.3.** According to the theorems 3.1, 3.2, mapping

$$f : H \to H \quad f = a \circ E + a \star E \quad a \in H$$
$$f \circ x = ax + xa$$

has matrix

$$\begin{pmatrix} a^0 & -a^1 & -a^2 & -a^3 \\ a^1 & a^0 & -a^3 & a^2 \\ a^2 & a^3 & a^0 & -a^1 \\ a^3 & -a^2 & a^1 & a^0 \end{pmatrix} + \begin{pmatrix} a^0 & -a^1 & -a^2 & -a^3 \\ a^1 & a^0 & a^3 & -a^2 \\ a^2 & -a^3 & a^0 & a^1 \\ a^3 & a^2 & -a^1 & a^0 \end{pmatrix}$$

$$= \begin{pmatrix} 2a^0 & -2a^1 & -2a^2 & -2a^3 \\ 2a^1 & 2a^0 & 0 & 0 \\ 2a^2 & 0 & 2a^0 & 0 \\ 2a^3 & 0 & 0 & 2a^0 \end{pmatrix}$$

According to the theorem 7.1,

$$a_0^0 = \frac{1}{2}(-2a^0 + 2a^0 + 2a^0 + 2a^0) \qquad a_1^0 = \frac{1}{2}(2a^0 - 2a^0)$$
$$a_0^1 = \frac{1}{2}(-2a^1 - (-2a^1) + 0 - 0) \qquad a_1^1 = \frac{1}{2}(2a^1 + (-2a^1))$$
$$a_0^2 = \frac{1}{2}(-2a^2 - 0 - (-2a^2) + 0) \qquad a_1^2 = \frac{1}{2}(2a^2 + 0)$$
$$a_0^3 = \frac{1}{2}(-2a^3 + 0 - 0 - (-2a^3)) \qquad a_1^3 = \frac{1}{2}(2a^3 - 0)$$

$$a_2^0 = \frac{1}{2}(2a^0 - 2a^0) \qquad a_3^0 = \frac{1}{2}(2a^0 - 2a^0)$$
$$a_2^1 = \frac{1}{2}(2a^1 - 0) \qquad a_3^1 = \frac{1}{2}(2a^1 + 0)$$
$$a_2^2 = \frac{1}{2}(2a^2 + (-2a^2)) \qquad a_3^2 = \frac{1}{2}(2a^2 - 0)$$
$$a_2^3 = \frac{1}{2}(2a^3 + 0) \qquad a_3^3 = \frac{1}{2}(2a^3 + (-2a^3))$$

Therefore

$$
\begin{array}{llll}
a_0^{\mathbf{0}} = 2a^{\mathbf{0}} & a_1^{\mathbf{0}} = 0 & a_2^{\mathbf{0}} = 0 & a_3^{\mathbf{0}} = 0 \\
a_0^{\mathbf{1}} = 0 & a_1^{\mathbf{1}} = 0 & a_2^{\mathbf{1}} = a^{\mathbf{1}} & a_3^{\mathbf{1}} = a^{\mathbf{1}} \\
a_0^{\mathbf{2}} = 0 & a_1^{\mathbf{2}} = a^{\mathbf{2}} & a_2^{\mathbf{2}} = 0 & a_3^{\mathbf{2}} = a^{\mathbf{2}} \\
a_0^{\mathbf{3}} = 0 & a_1^{\mathbf{3}} = a^{\mathbf{3}} & a_2^{\mathbf{3}} = a^{\mathbf{3}} & a_3^{\mathbf{3}} = 0
\end{array}
$$

□

## 9. Special Symbols and Notations

$E_{l \cdot a}$  Jacobian matrix of left shift 2

$\overline{E}$  linear automorphism of quaternioin algebra 2

$E_{r \cdot a}$  Jacobian matrix of right shift 3

# Отображения сопряжения алгебры кватернионов

Александр Клейн

Аннотация. В статье рассмотрены отображения сопряжения алгебры кватернионов. Доказана теорема о единственности разложения $R$-линейного отображения алгебры кватернионов относительно заданного множества отображений сопряжения.

Содержание

## 1. Предисловие

Теоремы [1]-2.5.1, [3]-1.5.1 утверждают, что $R$-линейное отображение поля комплексных чисел имеет вид

(1.1) $$f = a + b \circ I$$

(1.2) $$(a + b \circ I) \circ z = az + b\overline{z}$$

где

$$I \circ z = \overline{z}$$

сопряжение поля комплексных чисел. Тождественное отображение является линейным гомоморфизмом. Сопряжение является антилинейным гомоморфизмом.

В статье [2], я обнаружил, что существуют нетривиальные линейные автоморфизмы алгебры кватернионов. Это вдохновило меня предположить, что возможно разложить $R$-линейное отображение в сумму линейных и антилинейных автоморфизмов. Такое разложение существует (теорема [2]-12.1 ). Однако такое разложение линейного отображения не позволяет выделить в множестве

Aleks_Kleyn@MailAPS.org.
http://sites.google.com/site/AleksKleyn/.
http://arxiv.org/a/kleyn_a_1.
http://AleksKleyn.blogspot.com/.

линейных отображений класс отображений, которые были бы похожи на линейный гомоморфизм (пример [2]-12.3 ).

Размерность поля комплексных чисел равна 2, однако размерность алгебры кватернионов равна 4. Поэтому найти удовлетворительное разложение линейного отображения относительно базиса, состоящего из линейных и антилинейных гомоморфизмов, задача непростая.

Отображение $I$ в равенстве (1.1) является сопряжением. Однако мы не можем найти разложение линейного отображения в сумму тождественного отображения и сопряжения, так как этих отображений недостаточно для построения базиса.

Сопряжение является отображением симметрии. При сопряжении, коэффициенты при $i$, $j$, $k$ меняют знак. В этой статье я хочу изучить как будет выглядеть разложение линейного отображения, если я в качестве базиса выберу отображения

$$(\overline{E},\overline{I},\overline{J},\overline{K})$$

где каждое из рассматриваемых отображений меняет знак только одного коэффициента.

## 2. Соглашения

**Соглашение 2.1.** *Пусть $A$ - свободная конечно мерная алгебра. При разложении элемента алгебры $A$ относительно базиса $\overline{\overline{e}}$ мы пользуемся одной и той же корневой буквой для обозначения этого элемента и его координат. Однако в алгебре не принято использовать векторные обозначения. В выражении $a^2$ не ясно - это компонента разложения элемента $a$ относительно базиса или это операция возведения в степень. Для облегчения чтения текста мы будем индекс элемента алгебры выделять цветом. Например,*

$$a = a^i \overline{e}_i$$

□

**Соглашение 2.2.** *Если свободная конечномерная алгебра имеет единицу, то мы будем отождествлять вектор базиса $\overline{e}_0$ с единицей алгебры.* □

Без сомнения, у читателя могут быть вопросы, замечания, возражения. Я буду признателен любому отзыву.

## 3. Отображение $\overline{E}$

Рассмотрим отображение

$$\overline{E}: H \to H \quad \overline{E} \circ x = x$$

$$E^i_j = \delta^i_j$$

**Теорема 3.1.** *Мы можем отождествить отображение*

$$a \circ \overline{E} : H \to H \quad a \in H$$

и матрицу

$$(3.1) \quad E_{l \cdot a} = \begin{pmatrix} a^0 & -a^1 & -a^2 & -a^3 \\ a^1 & a^0 & -a^3 & a^2 \\ a^2 & a^3 & a^0 & -a^1 \\ a^3 & -a^2 & a^1 & a^0 \end{pmatrix}$$

*Доказательство.* Произведение кватернионов

$$a = a^0 + a^1 i + a^2 j + a^3 k$$

и

$$x = x^0 + x^1 i + x^2 j + x^3 k$$

имеет вид

$$ax = a^0 x^0 - a^1 x^1 - a^2 x^2 - a^3 x^3 + (a^0 x^1 + a^1 x^0 + a^2 x^3 - a^3 x^2)i$$
$$+ (a^0 x^2 + a^2 x^0 + a^3 x^1 - a^1 x^3)j + (a^0 x^3 + a^3 x^0 + a^1 x^2 - a^2 x^1)k$$

Следовательно, отображение $f_a(x) = ax$ имеет матрицу Якоби (3.1). □

**Теорема 3.2.** *Мы можем отождествить отображение*

$$a \star \overline{E} : H \to H \quad a \in H$$

*и матрицу*

$$(3.2) \quad E_{r \cdot a} = \begin{pmatrix} a^0 & -a^1 & -a^2 & -a^3 \\ a^1 & a^0 & a^3 & -a^2 \\ a^2 & -a^3 & a^0 & a^1 \\ a^3 & a^2 & -a^1 & a^0 \end{pmatrix}$$

*Доказательство.* Произведение кватернионов

$$x = x^0 + x^1 i + x^2 j + x^3 k$$

и

$$a = a^0 + a^1 i + a^2 j + a^3 k$$

имеет вид

$$xa = x^0 a^0 - x^1 a^1 - x^2 a^2 - x^3 a^3 + (x^0 a^1 + x^1 a^0 + x^2 a^3 - x^3 a^2)i$$
$$+ (x^0 a^2 + x^2 a^0 + x^3 a^1 - x^1 a^3)j + (x^0 a^3 + x^3 a^0 + x^1 a^2 - x^2 a^1)k$$

Следовательно, отображение $f_a(x) = xa$ имеет матрицу Якоби (3.2). □

## 4. Отображение $\overline{I}$

Рассмотрим отображение

$$\overline{I} : H \to H \quad \overline{I} \circ x = x^0 - x^1 i + x^2 j + x^3 k$$

$$I_0^0 = 1 \quad I_1^1 = -1 \quad I_2^2 = 1 \quad I_3^3 = 1$$

$$I = \begin{pmatrix} 1 & 0 & 0 & 0 \\ 0 & -1 & 0 & 0 \\ 0 & 0 & 1 & 0 \\ 0 & 0 & 0 & 1 \end{pmatrix}$$

**Теорема 4.1.** *Мы можем отождествить отображение*

$$a \circ \overline{I} : H \to H \quad a \in H$$

*и матрицу*

(4.1) $$I_{l \cdot a} = \begin{pmatrix} a^0 & a^1 & -a^2 & -a^3 \\ a^1 & -a^0 & -a^3 & a^2 \\ a^2 & -a^3 & a^0 & -a^1 \\ a^3 & a^2 & a^1 & a^0 \end{pmatrix}$$

(4.2) $$I_{l \cdot a} = E_{l \cdot a \circ} {}^\circ I$$

*Доказательство.* Произведение кватернионов

$$a = a^0 + a^1 i + a^2 j + a^3 k$$

и

$$I \circ x = x^0 - x^1 i + x^2 j + x^3 k$$

имеет вид

$$a \circ I \circ x = a^0 x^0 + a^1 x^1 - a^2 x^2 - a^3 x^3 + (-a^0 x^1 + a^1 x^0 + a^2 x^3 - a^3 x^2) i$$
$$+ (a^0 x^2 + a^2 x^0 - a^3 x^1 - a^1 x^3) j + (a^0 x^3 + a^3 x^0 + a^1 x^2 + a^2 x^1) k$$

Следовательно, отображение $f_a \circ x = a \circ I \circ x$ имеет матрицу Якоби (4.1). Равенство (4.2) следует из цепочки равенств

$$E_{l \cdot a \circ} {}^\circ I = \begin{pmatrix} a^0 & -a^1 & -a^2 & -a^3 \\ a^1 & a^0 & -a^3 & a^2 \\ a^2 & a^3 & a^0 & -a^1 \\ a^3 & -a^2 & a^1 & a^0 \end{pmatrix} {}_\circ{}^\circ \begin{pmatrix} 1 & 0 & 0 & 0 \\ 0 & -1 & 0 & 0 \\ 0 & 0 & 1 & 0 \\ 0 & 0 & 0 & 1 \end{pmatrix}$$

$$= \begin{pmatrix} a^0 & a^1 & -a^2 & -a^3 \\ a^1 & -a^0 & -a^3 & a^2 \\ a^2 & -a^3 & a^0 & -a^1 \\ a^3 & a^2 & a^1 & a^0 \end{pmatrix} = I_{l \cdot a}$$

$\square$

**Теорема 4.2.** *Мы можем отождествить отображение*

$$a \star \overline{I} : H \to H \quad a \in H$$

и матрицу

(4.3) $$I_{r\cdot a} = \begin{pmatrix} a^0 & a^1 & -a^2 & -a^3 \\ a^1 & -a^0 & a^3 & -a^2 \\ a^2 & a^3 & a^0 & a^1 \\ a^3 & -a^2 & -a^1 & a^0 \end{pmatrix}$$

(4.4) $$I_{r\cdot a} = E_{r\cdot a\circ}{}^\circ I$$

*Доказательство.* Произведение кватернионов
$$I \circ x = x^0 - x^1 i + x^2 j + x^3 k$$
и
$$a = a^0 + a^1 i + a^2 j + a^3 k$$
имеет вид
$$(I\circ x)a = x^0 a^0 + x^1 a^1 - x^2 a^2 - x^3 a^3 + (x^0 a^1 - x^1 a^0 + x^2 a^3 - x^3 a^2)i$$
$$+ (x^0 a^2 + x^2 a^0 + x^3 a^1 + x^1 a^3)j + (x^0 a^3 + x^3 a^0 - x^1 a^2 - x^2 a^1)k$$

Следовательно, отображение $f_a \circ x = a \star I \circ x$ имеет матрицу Якоби (4.3). Равенство (4.4) следует из цепочки равенств

$$E_{r\cdot a\circ}{}^\circ I = \begin{pmatrix} a^0 & -a^1 & -a^2 & -a^3 \\ a^1 & a^0 & a^3 & -a^2 \\ a^2 & -a^3 & a^0 & a^1 \\ a^3 & a^2 & -a^1 & a^0 \end{pmatrix} \circ^\circ \begin{pmatrix} 1 & 0 & 0 & 0 \\ 0 & -1 & 0 & 0 \\ 0 & 0 & 1 & 0 \\ 0 & 0 & 0 & 1 \end{pmatrix}$$

$$= \begin{pmatrix} a^0 & a^1 & -a^2 & -a^3 \\ a^1 & -a^0 & a^3 & -a^2 \\ a^2 & a^3 & a^0 & a^1 \\ a^3 & -a^2 & -a^1 & a^0 \end{pmatrix} = I_{r\cdot a}$$

□

## 5. Отображение $\overline{J}$

Рассмотрим отображение
$$\overline{J}: H \to H \quad \overline{J} \circ x = x^0 + x^1 i - x^2 j + x^3 k$$
$$J^0_0 = 1 \quad J^1_1 = 1 \quad J^2_2 = -1 \quad J^3_3 = 1$$
$$J = \begin{pmatrix} 1 & 0 & 0 & 0 \\ 0 & 1 & 0 & 0 \\ 0 & 0 & -1 & 0 \\ 0 & 0 & 0 & 1 \end{pmatrix}$$

**Теорема 5.1.** *Мы можем отождествить отображение*

$$a \circ \overline{J} : H \to H \quad a \in H$$

*и матрицу*

(5.1) $$J_{l \cdot a} = \begin{pmatrix} a^0 & -a^1 & a^2 & -a^3 \\ a^1 & a^0 & a^3 & a^2 \\ a^2 & a^3 & -a^0 & -a^1 \\ a^3 & -a^2 & -a^1 & a^0 \end{pmatrix}$$

(5.2) $$J_{l \cdot a} = E_{l \cdot a \circ}{}^\circ J$$

*Доказательство.* Произведение кватернионов

$$a = a^0 + a^1 i + a^2 j + a^3 k$$

и

$$J \circ x = x^0 + x^1 i - x^2 j + x^3 k$$

имеет вид

$$a \circ J \circ x = a^0 x^0 - a^1 x^1 + a^2 x^2 - a^3 x^3 + (a^0 x^1 + a^1 x^0 + a^2 x^3 + a^3 x^2)i$$
$$+ (-a^0 x^2 + a^2 x^0 + a^3 x^1 - a^1 x^3)j + (a^0 x^3 + a^3 x^0 - a^1 x^2 - a^2 x^1)k$$

Следовательно, отображение $f_a \circ x = a \circ J \circ x$ имеет матрицу Якоби (5.1). Равенство (5.2) следует из цепочки равенств

$$E_{l \cdot a \circ}{}^\circ J = \begin{pmatrix} a^0 & -a^1 & -a^2 & -a^3 \\ a^1 & a^0 & -a^3 & a^2 \\ a^2 & a^3 & a^0 & -a^1 \\ a^3 & -a^2 & a^1 & a^0 \end{pmatrix} \circ{}^\circ \begin{pmatrix} 1 & 0 & 0 & 0 \\ 0 & 1 & 0 & 0 \\ 0 & 0 & -1 & 0 \\ 0 & 0 & 0 & 1 \end{pmatrix}$$

$$= \begin{pmatrix} a^0 & -a^1 & a^2 & -a^3 \\ a^1 & a^0 & a^3 & a^2 \\ a^2 & a^3 & -a^0 & -a^1 \\ a^3 & -a^2 & -a^1 & a^0 \end{pmatrix} = J_{l \cdot a}$$

□

**Теорема 5.2.** *Мы можем отождествить отображение*

$$a \star \overline{J} : H \to H \quad a \in H$$

*и матрицу*

(5.3) $$J_{r \cdot a} = \begin{pmatrix} a^0 & -a^1 & a^2 & -a^3 \\ a^1 & a^0 & -a^3 & -a^2 \\ a^2 & -a^3 & -a^0 & a^1 \\ a^3 & a^2 & a^1 & a^0 \end{pmatrix}$$

(5.4) $$J_{r \cdot a} = E_{r \cdot a \circ} {}^\circ J$$

*Доказательство.* Произведение кватернионов

$$J \circ x = x^0 + x^1 i - x^2 j + x^3 k$$

и

$$a = a^0 + a^1 i + a^2 j + a^3 k$$

имеет вид

$$(J \circ x)a = x^0 a^0 - x^1 a^1 + x^2 a^2 - x^3 a^3 + (x^0 a^1 + x^1 a^0 - x^2 a^3 - x^3 a^2)i$$
$$+ (x^0 a^2 - x^2 a^0 + x^3 a^1 - x^1 a^3)j + (x^0 a^3 + x^3 a^0 + x^1 a^2 + x^2 a^1)k$$

Следовательно, отображение $f_a \circ x = a \star J \circ x$ имеет матрицу Якоби (5.3). Равенство (5.4) следует из цепочки равенств

$$E_{r \cdot a \circ} {}^\circ I = \begin{pmatrix} a^0 & -a^1 & -a^2 & -a^3 \\ a^1 & a^0 & a^3 & -a^2 \\ a^2 & -a^3 & a^0 & a^1 \\ a^3 & a^2 & -a^1 & a^0 \end{pmatrix} \circ^\circ \begin{pmatrix} 1 & 0 & 0 & 0 \\ 0 & -1 & 0 & 0 \\ 0 & 0 & 1 & 0 \\ 0 & 0 & 0 & 1 \end{pmatrix}$$

$$= \begin{pmatrix} a^0 & a^1 & -a^2 & -a^3 \\ a^1 & -a^0 & a^3 & -a^2 \\ a^2 & a^3 & a^0 & a^1 \\ a^3 & -a^2 & -a^1 & a^0 \end{pmatrix} = I_{r \cdot a}$$

□

## 6. Отображение $\overline{K}$

Рассмотрим отображение

$$\overline{K} : H \to H \quad \overline{K} \circ x = x^0 + x^1 i + x^2 j - x^3 k$$

$$K_0^0 = 1 \quad K_1^1 = 1 \quad K_2^2 = 1 \quad K_3^3 = -1$$

$$K = \begin{pmatrix} 1 & 0 & 0 & 0 \\ 0 & 1 & 0 & 0 \\ 0 & 0 & 1 & 0 \\ 0 & 0 & 0 & -1 \end{pmatrix}$$

**Теорема 6.1.** *Мы можем отождествить отображение*

$$a \circ \overline{K} : H \to H \quad a \in H$$

и матрицу

$$K_{l \cdot a} = \begin{pmatrix} a^0 & -a^1 & -a^2 & a^3 \\ a^1 & a^0 & -a^3 & -a^2 \\ a^2 & a^3 & a^0 & a^1 \\ a^3 & -a^2 & a^1 & -a^0 \end{pmatrix} \tag{6.1}$$

$$K_{l \cdot a} = E_{l \cdot a \circ}{}^{\circ} K \tag{6.2}$$

*Доказательство.* Произведение кватернионов

$$a = a^0 + a^1 i + a^2 j + a^3 k$$

и

$$K \circ x = x^0 + x^1 i + x^2 j - x^3 k$$

имеет вид

$$a \circ K \circ x = a^0 x^0 - a^1 x^1 - a^2 x^2 + a^3 x^3 + (a^0 x^1 + a^1 x^0 - a^2 x^3 - a^3 x^2)i$$
$$+ (a^0 x^2 + a^2 x^0 + a^3 x^1 + a^1 x^3)j + (-a^0 x^3 + a^3 x^0 + a^1 x^2 - a^2 x^1)k$$

Следовательно, отображение $f_a \circ x = a \circ K \circ x$ имеет матрицу Якоби (6.1). Равенство (6.2) следует из цепочки равенств

$$E_{l \cdot a \circ}{}^{\circ} K = \begin{pmatrix} a^0 & -a^1 & -a^2 & -a^3 \\ a^1 & a^0 & -a^3 & a^2 \\ a^2 & a^3 & a^0 & -a^1 \\ a^3 & -a^2 & a^1 & a^0 \end{pmatrix} {}_{\circ}{}^{\circ} \begin{pmatrix} 1 & 0 & 0 & 0 \\ 0 & 1 & 0 & 0 \\ 0 & 0 & 1 & 0 \\ 0 & 0 & 0 & -1 \end{pmatrix}$$

$$= \begin{pmatrix} a^0 & -a^1 & -a^2 & a^3 \\ a^1 & a^0 & -a^3 & -a^2 \\ a^2 & a^3 & a^0 & a^1 \\ a^3 & -a^2 & a^1 & -a^0 \end{pmatrix} = K_{l \cdot a}$$

□

**Теорема 6.2.** *Мы можем отождествить отображение*

$$a \star \overline{K} : H \to H \quad a \in H$$

*и матрицу*

$$K_{r \cdot a} = \begin{pmatrix} a^0 & -a^1 & -a^2 & a^3 \\ a^1 & a^0 & a^3 & a^2 \\ a^2 & -a^3 & a^0 & -a^1 \\ a^3 & a^2 & -a^1 & -a^0 \end{pmatrix} \tag{6.3}$$

$$K_{r \cdot a} = E_{r \cdot a \circ}{}^{\circ} K \tag{6.4}$$

*Доказательство.* Произведение кватернионов
$$K \circ x = x^0 + x^1 i + x^2 j - x^3 k$$
и
$$a = a^0 + a^1 i + a^2 j + a^3 k$$
имеет вид
$$(K \circ x)a = x^0 a^0 - x^1 a^1 - x^2 a^2 + x^3 a^3 + (x^0 a^1 + x^1 a^0 + x^2 a^3 + x^3 a^2)i$$
$$+ (x^0 a^2 + x^2 a^0 - x^3 a^1 - x^1 a^3)j + (x^0 a^3 - x^3 a^0 + x^1 a^2 - x^2 a^1)k$$

Следовательно, отображение $f_a \circ x = a \star K \circ x$ имеет матрицу Якоби (6.3). Равенство (6.4) следует из цепочки равенств

$$E_{r \cdot a \circ} \circ I = \begin{pmatrix} a^0 & -a^1 & -a^2 & -a^3 \\ a^1 & a^0 & a^3 & -a^2 \\ a^2 & -a^3 & a^0 & a^1 \\ a^3 & a^2 & -a^1 & a^0 \end{pmatrix} \circ \begin{pmatrix} 1 & 0 & 0 & 0 \\ 0 & -1 & 0 & 0 \\ 0 & 0 & 1 & 0 \\ 0 & 0 & 0 & 1 \end{pmatrix}$$

$$= \begin{pmatrix} a^0 & a^1 & -a^2 & -a^3 \\ a^1 & -a^0 & a^3 & -a^2 \\ a^2 & a^3 & a^0 & a^1 \\ a^3 & -a^2 & -a^1 & a^0 \end{pmatrix} = I_{r \cdot a}$$

□

## 7. Линейное отображение

**Теорема 7.1.** *Линейное отображение алгебры кватернионов*
$$\overline{f} : H \to H$$
*имеет единственное разложение*

(7.1) $$\overline{f} = a_0 \circ \overline{E} + a_1 \circ \overline{I} + a_2 \circ \overline{J} + a_3 \circ \overline{K}$$

*При этом*

(7.2) $$a_0^0 = \frac{1}{2}(-f_0^0 + f_1^1 + f_2^2 + f_3^3)$$

(7.3) $$a_1^0 = \frac{1}{2}(f_0^0 - f_1^1)$$

(7.4) $$a_2^0 = \frac{1}{2}(f_0^0 - f_2^2)$$

(7.5) $$a_3^0 = \frac{1}{2}(f_0^0 - f_3^3)$$

$$(7.6) \quad a_0^1 = \frac{1}{2}(-f_0^1 - f_1^0 + f_2^3 - f_3^2)$$

$$(7.7) \quad a_1^1 = \frac{1}{2}(f_0^1 + f_1^0)$$

$$(7.8) \quad a_2^1 = \frac{1}{2}(f_0^1 - f_2^3)$$

$$(7.9) \quad a_3^1 = \frac{1}{2}(f_0^1 + f_3^2)$$

$$(7.10) \quad a_0^2 = \frac{1}{2}(-f_0^2 - f_1^3 - f_2^0 + f_3^1)$$

$$(7.11) \quad a_1^2 = \frac{1}{2}(f_0^2 + f_1^3)$$

$$(7.12) \quad a_2^2 = \frac{1}{2}(f_0^2 + f_2^0)$$

$$(7.13) \quad a_3^2 = \frac{1}{2}(f_0^2 - f_3^1)$$

$$(7.14) \quad a_0^3 = \frac{1}{2}(-f_0^3 + f_1^2 - f_2^1 - f_3^0)$$

$$(7.15) \quad a_1^3 = \frac{1}{2}(f_0^3 - f_1^2)$$

$$(7.16) \quad a_2^3 = \frac{1}{2}(f_0^3 + f_2^1)$$

$$(7.17) \quad a_3^3 = \frac{1}{2}(f_0^3 + f_3^0)$$

*Доказательство.* Линейное отображение (7.1) имеет матрицу

$$(7.18) \quad \begin{pmatrix} a_0^0 & -a_0^1 & -a_0^2 & -a_0^3 \\ a_0^1 & a_0^0 & -a_0^3 & a_0^2 \\ a_0^2 & a_0^3 & a_0^0 & -a_0^1 \\ a_0^3 & -a_0^2 & a_0^1 & a_0^0 \end{pmatrix} + \begin{pmatrix} a_1^0 & a_1^1 & -a_1^2 & -a_1^3 \\ a_1^1 & -a_1^0 & -a_1^3 & a_1^2 \\ a_1^2 & -a_1^3 & a_1^0 & -a_1^1 \\ a_1^3 & a_1^2 & a_1^1 & a_1^0 \end{pmatrix} \\ + \begin{pmatrix} a_2^0 & -a_2^1 & a_2^2 & -a_2^3 \\ a_2^1 & a_2^0 & a_2^3 & a_2^2 \\ a_2^2 & a_2^3 & -a_2^0 & -a_2^1 \\ a_2^3 & -a_2^2 & -a_2^1 & a_2^0 \end{pmatrix} + \begin{pmatrix} a_3^0 & -a_3^1 & -a_3^2 & a_3^3 \\ a_3^1 & a_3^0 & -a_3^3 & -a_3^2 \\ a_3^2 & a_3^3 & a_3^0 & a_3^1 \\ a_3^3 & -a_3^2 & a_3^1 & -a_3^0 \end{pmatrix}$$

Из сравнения матрицы отображения $\overline{f}$ и матрицы (7.18), мы получим систему линейных уравнений

$$(7.19) \quad f^0_0 = a^0_0 + a^0_1 + a^0_2 + a^0_3$$
$$(7.20) \quad f^0_1 = -a^1_0 + a^1_1 - a^1_2 - a^1_3$$
$$(7.21) \quad f^0_2 = -a^2_0 - a^2_1 + a^2_2 - a^2_3$$
$$(7.22) \quad f^0_3 = -a^3_0 - a^3_1 - a^3_2 + a^3_3$$
$$(7.23) \quad f^1_0 = a^1_0 + a^1_1 + a^1_2 + a^1_3$$
$$(7.24) \quad f^1_1 = a^0_0 - a^0_1 + a^0_2 + a^0_3$$
$$(7.25) \quad f^1_2 = -a^3_0 - a^3_1 + a^3_2 - a^3_3$$
$$(7.26) \quad f^1_3 = a^2_0 + a^2_1 + a^2_2 - a^2_3$$
$$(7.27) \quad f^2_0 = a^2_0 + a^2_1 + a^2_2 + a^2_3$$
$$(7.28) \quad f^2_1 = a^3_0 - a^3_1 + a^3_2 + a^3_3$$
$$(7.29) \quad f^2_2 = a^0_0 + a^0_1 - a^0_2 + a^0_3$$
$$(7.30) \quad f^2_3 = -a^1_0 - a^1_1 - a^1_2 + a^1_3$$
$$(7.31) \quad f^3_0 = a^3_0 + a^3_1 + a^3_2 + a^3_3$$
$$(7.32) \quad f^3_1 = -a^2_0 + a^2_1 - a^2_2 - a^2_3$$
$$(7.33) \quad f^3_2 = a^1_0 + a^1_1 - a^1_2 + a^1_3$$
$$(7.34) \quad f^3_3 = a^0_0 + a^0_1 + a^0_2 - a^0_3$$

Уравнение (7.3) следует из уравнений (7.19), (7.24). Уравнение (7.4) следует из уравнений (7.19), (7.29). Уравнение (7.5) следует из уравнений (7.19), (7.34). Уравнение (7.2) следует из уравнений (7.3), (7.4), (7.5), (7.19).

Уравнение (7.7) следует из уравнений (7.23), (7.20). Уравнение (7.8) следует из уравнений (7.23), (7.33). Уравнение (7.9) следует из уравнений (7.23), (7.30). Уравнение (7.6) следует из уравнений (7.7), (7.8), (7.9), (7.23).

Уравнение (7.11) следует из уравнений (7.27), (7.32). Уравнение (7.12) следует из уравнений (7.27), (7.21). Уравнение (7.13) следует из уравнений (7.27), (7.26). Уравнение (7.10) следует из уравнений (7.11), (7.12), (7.13), (7.27).

Уравнение (7.15) следует из уравнений (7.31), (7.28). Уравнение (7.16) следует из уравнений (7.31), (7.25). Уравнение (7.17) следует из уравнений (7.31), (7.22). Уравнение (7.14) следует из уравнений (7.15), (7.16), (7.17), (7.31). □

Из теоремы [3]-2.6.4 (смотри также теорему [1]-3.7.4) следует, что множество линейных эндоморфизмов $\mathcal{L}(H;H)$ алгебры кватернионов $H$ изоморфно тензорному произведению $H \otimes H$. Теорема 7.1 утверждает, что мы можем рассматривать модуль $\mathcal{L}(H;H)$ как $H\star$-векторное пространство с базисом

$$(\overline{E}, \overline{I}, \overline{J}, \overline{K})$$

**Пример 7.2.** Согласно теореме 3.2, отображение

$$a \star \overline{E} : H \to H \quad a \in H$$

имеет матрицу

$$\begin{pmatrix} a^0 & -a^1 & -a^2 & -a^3 \\ a^1 & a^0 & a^3 & -a^2 \\ a^2 & -a^3 & a^0 & a^1 \\ a^3 & a^2 & -a^1 & a^0 \end{pmatrix}$$

Согласно теореме 7.1,

$$a_0^0 = \frac{1}{2}(-a^0 + a^0 + a^0 + a^0) \qquad a_1^0 = \frac{1}{2}(a^0 - a^0)$$

$$a_0^1 = \frac{1}{2}(-a^1 - (-a^1) + (-a^1) - a^1) \qquad a_1^1 = \frac{1}{2}(a^1 + (-a^1))$$

$$a_0^2 = \frac{1}{2}(-a^2 - a^2 - (-a^2) + (-a^2)) \qquad a_1^2 = \frac{1}{2}(a^2 + a^2)$$

$$a_0^3 = \frac{1}{2}(-a^3 + (-a^3) - a^3 - (-a^3)) \qquad a_1^3 = \frac{1}{2}(a^3 - (-a^3))$$

$$a_2^0 = \frac{1}{2}(a^0 - a^0) \qquad a_3^0 = \frac{1}{2}(a^0 - a^0)$$

$$a_2^1 = \frac{1}{2}(a^1 - (-a^1)) \qquad a_3^1 = \frac{1}{2}(a^1 + a^1)$$

$$a_2^2 = \frac{1}{2}(a^2 + (-a^2)) \qquad a_3^2 = \frac{1}{2}(a^2 - (-a^2))$$

$$a_2^3 = \frac{1}{2}(a^3 + a^3) \qquad a_3^3 = \frac{1}{2}(a^3 + (-a^3))$$

Следовательно

$$\begin{array}{llll} a_0^0 = a^0 & a_1^0 = 0 & a_2^0 = 0 & a_3^0 = 0 \\ a_0^1 = -a^1 & a_1^1 = 0 & a_2^1 = a^1 & a_3^1 = a^1 \\ a_0^2 = -a^2 & a_1^2 = a^2 & a_2^2 = 0 & a_3^2 = a^2 \\ a_0^3 = -a^3 & a_1^3 = a^3 & a_2^3 = a^3 & a_3^3 = 0 \end{array}$$

$\square$

**Пример 7.3.** Согласно теоремам 3.1, 3.2, отображение

$$f : H \to H \quad f = a \circ E + a \star E \quad a \in H$$
$$f \circ x = ax + xa$$

имеет матрицу

$$\begin{pmatrix} a^0 & -a^1 & -a^2 & -a^3 \\ a^1 & a^0 & -a^3 & a^2 \\ a^2 & a^3 & a^0 & -a^1 \\ a^3 & -a^2 & a^1 & a^0 \end{pmatrix} + \begin{pmatrix} a^0 & -a^1 & -a^2 & -a^3 \\ a^1 & a^0 & a^3 & -a^2 \\ a^2 & -a^3 & a^0 & a^1 \\ a^3 & a^2 & -a^1 & a^0 \end{pmatrix}$$

$$= \begin{pmatrix} 2a^0 & -2a^1 & -2a^2 & -2a^3 \\ 2a^1 & 2a^0 & 0 & 0 \\ 2a^2 & 0 & 2a^0 & 0 \\ 2a^3 & 0 & 0 & 2a^0 \end{pmatrix}$$

Согласно теореме 7.1,

$$a_0^{\mathbf{0}} = \frac{1}{2}(-2a^{\mathbf{0}} + 2a^{\mathbf{0}} + 2a^{\mathbf{0}} + 2a^{\mathbf{0}}) \qquad a_1^{\mathbf{0}} = \frac{1}{2}(2a^{\mathbf{0}} - 2a^{\mathbf{0}})$$

$$a_0^{\mathbf{1}} = \frac{1}{2}(-2a^{\mathbf{1}} - (-2a^{\mathbf{1}}) + 0 - 0) \qquad a_1^{\mathbf{1}} = \frac{1}{2}(2a^{\mathbf{1}} + (-2a^{\mathbf{1}}))$$

$$a_0^{\mathbf{2}} = \frac{1}{2}(-2a^{\mathbf{2}} - 0 - (-2a^{\mathbf{2}}) + 0) \qquad a_1^{\mathbf{2}} = \frac{1}{2}(2a^{\mathbf{2}} + 0)$$

$$a_0^{\mathbf{3}} = \frac{1}{2}(-2a^{\mathbf{3}} + 0 - 0 - (-2a^{\mathbf{3}})) \qquad a_1^{\mathbf{3}} = \frac{1}{2}(2a^{\mathbf{3}} - 0)$$

$$a_2^{\mathbf{0}} = \frac{1}{2}(2a^{\mathbf{0}} - 2a^{\mathbf{0}}) \qquad a_3^{\mathbf{0}} = \frac{1}{2}(2a^{\mathbf{0}} - 2a^{\mathbf{0}})$$

$$a_2^{\mathbf{1}} = \frac{1}{2}(2a^{\mathbf{1}} - 0) \qquad a_3^{\mathbf{1}} = \frac{1}{2}(2a^{\mathbf{1}} + 0)$$

$$a_2^{\mathbf{2}} = \frac{1}{2}(2a^{\mathbf{2}} + (-2a^{\mathbf{2}})) \qquad a_3^{\mathbf{2}} = \frac{1}{2}(2a^{\mathbf{2}} - 0)$$

$$a_2^{\mathbf{3}} = \frac{1}{2}(2a^{\mathbf{3}} + 0) \qquad a_3^{\mathbf{3}} = \frac{1}{2}(2a^{\mathbf{3}} + (-2a^{\mathbf{3}}))$$

Следовательно

$$a_0^{\mathbf{0}} = 2a^{\mathbf{0}} \qquad a_1^{\mathbf{0}} = 0 \qquad a_2^{\mathbf{0}} = 0 \qquad a_3^{\mathbf{0}} = 0$$
$$a_0^{\mathbf{1}} = 0 \qquad a_1^{\mathbf{1}} = 0 \qquad a_2^{\mathbf{1}} = a^{\mathbf{1}} \qquad a_3^{\mathbf{1}} = a^{\mathbf{1}}$$
$$a_0^{\mathbf{2}} = 0 \qquad a_1^{\mathbf{2}} = a^{\mathbf{2}} \qquad a_2^{\mathbf{2}} = 0 \qquad a_3^{\mathbf{2}} = a^{\mathbf{2}}$$
$$a_0^{\mathbf{3}} = 0 \qquad a_1^{\mathbf{3}} = a^{\mathbf{3}} \qquad a_2^{\mathbf{3}} = a^{\mathbf{3}} \qquad a_3^{\mathbf{3}} = 0$$

□

## 8. Список литературы

# 9. Специальные символы и обозначения

$E_{l \cdot a}$  матрица Якоби левого сдвига 3

$\overline{E}$  линейный автоморфизм алгебры
   кватернионов 2

$E_{r \cdot a}$  матрица Якоби правого сдвига 3